\documentclass[11pt]{article}
\usepackage{amsmath,amssymb}

\def\r{\rightarrow}

\def\qed{\hfill\vrule height5pt width5pt depth0pt}
\def\one #1{1_{\{#1\}}}

\def\1{\mbox{\bf 1}}

\newcommand{\proof}{\noindent {\bf Proof:\ }}

\newtheorem{Theorem}{Theorem}
\newtheorem{Lemma}{Lemma}

\newtheorem{Corollary}{Corollary}
\newtheorem{Remark}{Remark}

\begin{document}
\title{Decomposition results for stochastic storage processes and queues with alternating L\'evy inputs
}
\author{
Onno Boxma\thanks{EURANDOM and Department of Mathematics and
Computer Science; Eindhoven University
of Technology; P.O. Box 513; 5600 MB Eindhoven; The Netherlands ({\tt boxma@win.tue.nl})}
\and
Offer Kella\thanks{Department of Statistics; The Hebrew University of Jerusalem; Mount Scopus,
Jerusalem 91905; Israel ({\tt offer.kella@huji.ac.il})}
\thanks{Supported in part by grant 434/09 from the Israel Science
Foundation, the Vigevani Chair in Statistics and visitor grant No. 040.11.257 from The Netherlands Organisation for Scientific Research.}
}
\date{October 8, 2012}
\maketitle

\abstract{
In this paper we generalize known workload decomposition results for L\'evy queues with secondary jump inputs and queues with server vacations or service interruptions. Special cases are polling systems with either compound Poisson or more general L\'evy inputs.
Our main tools are new martingale results, which have been derived in a companion paper.
}

\vspace{0.1in}

\vspace{0.1in}


\section{Introduction}
It is
well known in queueing theory (e.g., \cite{fc1985,s1988}) that in a stable M/G/1 queue with server down periods (vacations, interruptions, etc.) the steady state waiting time distribution (properly defined) is a convolution of two or more distributions, one of which is always the steady state waiting time distribution of an ordinary M/G/1 queue. As Poisson arrivals see time averages, this result also holds for the workload process.
\cite{kw1991} studies a more general model of a L\'evy process with no negative jumps and additional jumps that occur at stopping epochs and the size of which is measurable with respect to the current information.
The interesting outcome of \cite{kw1991} was that the same (and even more general) decomposition results that were known for queues also turned out to hold for these L\'evy processes with additional jumps.

That model is interesting in its own right but can also be viewed as a weak limit of queues with down (off) periods
where during these down periods workload can only accumulate as the server is idle.
Consider a process that can be either in an {\em up} (on) state or a {\em down} state. When it is in the {\em up} state it behaves like some L\'evy process with no negative jumps and a negative drift. When it is in a {\em down}
state it behaves like a subordinator, that is, a nondecreasing L\'evy process.
The question that comes to mind is whether this up/down process (for which we give a precise definition later) obeys a similar decomposition property. This would immediately imply a decomposition in certain polling systems as described in Section~\ref{sec:polling} below. It is a simple observation that if one cuts and pastes the up/down process such that only the up periods are visible, then the resulting process is the one that was considered in \cite{kw1991}. As it seemed that the results of \cite{kw1991} could not be used in our setting, we found it necessary to develop a more general theory, in particular a certain martingale theory that would streamline our work and could be useful in other applications as well.
That direction was developed in \cite{KB}.


 The first main result of \cite{KB} is the extension of the martingale results of \cite{kw1992} to the case where the driving process is a L\'evy-type process. That is, it is a sum of stochastic integrals of some bounded left continuous right limit process with respect to coordinate processes associated with some multidimensional L\'evy process.
Such processes with an even more general (predictable) integrand are discussed in \cite{applebaum}.
The second main result of \cite{KB} is that our local martingale is in fact an $L^2$ martingale, and moreover, when upon dividing by the time parameter $t$ it converges to zero almost surely and in $L^2$ as $t\r\infty$.

The main goal of the present paper is to apply the martingale results
which were derived in \cite{KB}
to establish decomposition results for the up/down model
that was introduced above.

The paper is organized as follows. In Theorem~\ref{martingale} of  Section~\ref{sec:prelim} we summarize the main results
from \cite{KB} which are needed in the present paper.
In Section~\ref{sec:4} we apply our results to establish decomposition results for the up/down model, thereby considerably generalizing the results of \cite{kw1991}.
In Section~\ref{sec:Wd} we identify the non-standard component in the decomposition associated with down periods. Finally in Section~\ref{sec:polling} a discussion of polling systems, the motivation for this study, is given and the contribution of our results to this area is emphasized.

For extensive discussions of decomposition results in queues and storage processes,
we refer to the surveys \cite{boxma,doshi} and to the recent study \cite{IK} and references therein.

\section{Preliminaries}
\label{sec:prelim}
In preparation of our analysis and in order to make this paper more self contained we first recall the main results from \cite{KB} which we will need here.

For what follows, given a c\`adl\`ag (right continuous left limit) function $g:\mathbb{R}_+\r\mathbb{R}$, we denote $g(t-)=\displaystyle \lim_{s\uparrow t}g(s)$, $\Delta g(t)=g(t)-g(t-)$ with the convention that $\Delta g(0)=g(0)$ and if $g$ is VF (finite variation on finite intervals), then $g^d(t)=\sum_{0\le s\le t}\Delta g(s)$ and $g^c(t)=g(t)-g^d(t)$. Also, $\mathbb{R}_+=[0,\infty)$, $\mathbb{R} =(-\infty,\infty)$ and {\em a.s.} abbreviates {\em almost surely}. Finally $a\wedge b=\min(a,b)$, $a\vee b=\max(a,b)$, $a^-=a\wedge 0$ and $a^+=a\vee0$.

Let $X=(X_1,\ldots,X_K)$ be a c\`adl\`ag $K$-dimensional L\'evy process with respect to some standard filtration $\left\{\mathcal{F}_t|\ t\ge0\right\}$ having no negative jumps (the L\'evy measure is concentrated on $\mathbb{R}_+^K$) with Laplace-Stieltjes exponent
\begin{eqnarray}\label{LSE}
\psi(\gamma)&=&\log Ee^{-\gamma^TX(1)}=-c^T\gamma+\frac{\gamma^T\Sigma\gamma}{2}\\ &&+\int_{\mathbb{R}_+^K}\left(e^{-\gamma^Tx}-1+\gamma^Tx\one{\|x\|\le 1}\right)\nu({\rm d}x)\ ,\nonumber
\end{eqnarray}
where $\gamma\in\mathbb{R}_+^K$. It is well known that in this case $\psi(\gamma)$ is finite for each $\gamma\ge 0$, convex (thus continuous) with $\psi(0)=0$, infinitely differentiable in the interior of $\mathbb{R}_+^K$, and for every $\gamma\ge 0$ for which $\gamma^TX$ is not a subordinator (not nondecreasing), $\psi(t\gamma)\r\infty$ as $t\r\infty$. Furthermore, $EX_k(t)=-t\frac{\partial\psi}{\partial\gamma_k}(0+)$ (finite or $+\infty$, but can never be $-\infty$) and when the first two right derivatives at zero are finite, then $\mbox{Cov}(X_k(t),X_\ell(t))=t\frac{\partial^2\psi}{\partial\gamma_k\partial\gamma_\ell}\psi(0+)$.

Let $I=(I_1,\ldots,I_k)$ be a nonnegative, bounded, c\`adl\`ag and adapted process and define (a special case of) a L\'evy-type process as a sum of the following stochastic integrals.
\begin{equation}
\tilde X(t)=\sum_{k=1}^K\int_{[0,t]}I_k(s-){\rm d}X_k(s).
\end{equation}
Finally, let $Y$ be a c\`adl\`ag adapted process having a.s. finite variation on finite intervals, set
\begin{equation}
Z(t)=\tilde X(t)+Y(t),
\end{equation}
and assume that $Z$ is bounded below. Under the above setup, the following summarizes what we need from~\cite{KB}.

\begin{Theorem}\label{martingale}
Given the assumption above:
\begin{description}
\item{(i)} \label{i} The following is a mean square martingale having zero mean:
\begin{eqnarray}\label{eq:LSTmart}
M(t)&=&\int_0^t\psi( I(s))e^{- Z(s)}{\rm d}s+e^{- Z(0)}-e^{- Z(t)}-\int_0^te^{- Z(s)}{\rm d}Y^c(s)\nonumber\\ \\
&&+\sum_{0<s\le t}e^{- Z(s)}\left(1-e^{ \Delta Y(s)}\right). \nonumber
\end{eqnarray}
\item{(ii)} \label{ii} $M(t)/t\to0$ as $t\to\infty$ a.s. and in $L^2$.
\item{(iii)} \label{iii} If
\begin{equation}
\frac{1}{t}\int_0^tI_k(s){\rm d}s\to \beta_k\ ,
\end{equation}
a.s. as $t\to\infty$ for each $k$, and if $EX_k(1)<\infty$ ($EX_k(1)^->-\infty$ as there are no negative jumps), then
\begin{equation}\label{eq:6}
\frac{\tilde X(t)}{t}\to \sum_{k=1}^K\beta_kEX_k(1)
\end{equation}
a.s. as $t\to\infty$.
\item{(iv)} \label{iv} When $\displaystyle Y(t)=-\inf_{0\le s\le t}\tilde X(s)^-$ and (\ref{eq:6}) holds, then
\begin{align}
\frac{1}{t}\int_0^t\psi(I(s))e^{-Z(s)}{\rm d}s\to-\left(\sum_{k=1}^K\beta_kEX_k(1)\right)^-\ ,
\end{align}
a.s. as $t\to\infty$.
\end{description}
\end{Theorem}

\section{Decomposition results for L\'evy storage processes}
\label{sec:4}

In this section we complement the results of \cite{kw1991} as follows. Let $0=T_0\le  S_1\le T_1\le S_2\le T_2\ldots$ be an increasing sequence of a.s. finite stopping times with respect to the standard filtration $\left\{\mathcal{F}_t|\ t\ge 0\right\}$ satisfying $T_{n-1}<T_n$ and $T_n\r\infty$ a.s. Let $X_n=S_n-T_{n-1}$ and $Y_n=T_n-S_n$. The model here is that $(T_{n-1},S_n]$ with lengths $X_n$ are {\em down} periods, where there is no output (the ``server'' is not working) and therefore the buffer content can only accumulate. $(S_n,T_n]$ with length $Y_n$ are {\em up} periods where there is both input and output, which is modeled as usual by a reflected (Skorohod map of the) process.
\begin{Remark}\rm
We note that in some models it is possible that there is no reflection. This occurs, for example, whenever the server is shut off as soon as the system empties, which may be modeled via the stopping times.
\end{Remark}
Let $X_u$ be a one-dimensional c\`adl\`ag L\'evy process with no negative jumps which is not a subordinator (not nondecreasing), and with Laplace-Stieltjes exponent
\begin{equation}
\varphi(\alpha)=-c_u\alpha+\frac{\sigma_u^2\alpha^2}{2}+\int_{(0,\infty)}\left(e^{-\alpha x}-1+\alpha x\one{x\le 1}\right)\nu_u({\rm d}x)
\end{equation}
 and assume that $EX_u(1)=-\varphi'(0)=c_u+\int_{(1,\infty)}x\nu_u({\rm d}x)<0$. This models the net input process (input minus potential output) during up periods. Let $X_d$ be a one-dimensional right continuous subordinator (nondecreasing L\'evy process) with Laplace-Stieltjes exponent $-\eta$ where
\begin{equation}
\eta(\alpha)=c_d\alpha+\int_{(0,\infty]}\left(1-e^{-\alpha x}\right)\nu_d({\rm d}x)
\end{equation}
and assume that $EX_d(1)=\eta'(0)<\infty$. The latter  models the process according to which work accumulates during down periods.

Now, set $N(t)=\sup\{n|\ T_n\le t\}$ and let $J(t)=\one{S_{N(t)+1}>t}$ and thus $J(t)=\one{J(t)=1}$ and $1-J(t)=\one{J(t)=0}$. That is, $J(t)=1$ during down periods and $J(t)=0$ during up periods. Finally, for $W(0)\in\mathcal{F}_0$ let
\begin{eqnarray}
\tilde X_d(t)&=&\int_{(0,t]}J(s-){\rm d}X_d(s),\nonumber\\
\tilde X_u(t)&=&\int_{(0,t]}(1-J(s-)){\rm d}X_u(s),\nonumber\\
\tilde X(t)&=&\tilde X_u(t)+\tilde X_d(t),\\
L(t)&=&\displaystyle-\inf_{0\le s\le t}(W(0)+\tilde X(s))^-,\nonumber\\
W(t)&=&W(0)+\tilde X(t)+L(t)\ .\nonumber
\end{eqnarray}
The process $\{W(t)|\ t\ge0\}$ is the content process of interest for which we would like to establish a general decomposition. During down periods it behaves like a subordinator with exponent $-\eta$ (and only grows) and during up periods it behaves like L\'evy process with exponent $\varphi$ and is reflected at the origin. This general decomposition is given by the following theorem which will be interpreted after its proof.

\begin{Theorem}\label{th:decomp} If, in addition to the above setup and assumptions,
\begin{equation}\label{eq:Winfty}
\frac{1}{t}\int_0^t e^{-\alpha W(s)}{\rm d}s\r Ee^{-\alpha W(\infty)}
\end{equation}
 a.s. as $t\r\infty$ (ergodic convergence) for some finite random variable $W(\infty)$ and
\begin{equation}\label{eq:jbar}
\frac{1}{t}\int_0^tJ(s){\rm d}s\r p_d\le \frac{\varphi'(0)}{\eta'(0)+\varphi'(0)}\ ,
\end{equation}
then there exists a nonnegative random variable $W_d$ such that if $p_d>0$
then a.s.
\begin{equation}
\frac{\int_0^t e^{-\alpha W(s)}J(s){\rm d}s}{\int_0^tJ(s){\rm d}s}\r Ee^{-\alpha W_d}
\end{equation}
for every $\alpha\ge 0$. Moreover, with
\begin{equation}
\pi_\ell=1-\left(1+\frac{\eta'(0)}{\varphi'(0)}\right) p_d
\end{equation}
and
\begin{equation}
\pi=\frac{\eta'(0)}{\eta'(0)+\varphi'(0)}
\end{equation}
we have that
\begin{align}\label{eq:decomp}
Ee^{-\alpha W(\infty)}&=\pi_\ell\frac{\varphi'(0)\alpha}{\varphi(\alpha)}\\ &\quad\nonumber+(1-\pi_\ell)\left(1-\pi+\pi\frac{\eta(\alpha)}{\eta'(0)\alpha}\frac{\varphi'(0)\alpha}{\varphi(\alpha)}\right) Ee^{-\alpha W_d} .
\end{align}
\end{Theorem}

\proof
With $\psi(\gamma_1,\gamma_2)=\varphi(\gamma_1)-\eta(\gamma_2)$ (or any other $\psi$ with $\psi(\alpha,0)=\varphi(\alpha)$ and $\psi(0,\alpha)=-\eta(\alpha)$), $I_1(s)=\alpha (1-J(s))$ and $I_2(s)=\alpha J(s)$,
Theorem~\ref{martingale}-(iv) and (\ref{eq:jbar}) imply that
\begin{align}\label{eq:converge}
&\frac{1}{t}\int_0^t(\varphi(\alpha)(1-J(s))-\eta(\alpha)J(s))e^{-\alpha W(s)}{\rm d}s\\
&=\varphi(\alpha)\frac{1}{t}\int_0^te^{-\alpha W(s)}{\rm d}s-(\varphi(\alpha)+\eta(\alpha))\frac{1}{t}\int_0^tJ(s)e^{-\alpha W(s)}{\rm d}s\nonumber
\end{align}
converges a.s., as $t\to\infty$, to
\begin{align}\label{eq:drift}
-\alpha(-(1-p_d)\varphi'(0)+p_d\eta'(0))^-=\alpha((1-p_d)\varphi'(0)-p_d\eta'(0))\ ,
\end{align}
the last equality being due to  $p_d\le \frac{\varphi'(0)}{\eta'(0)+\varphi'(0)}$.

Now, by the convergence of (\ref{eq:converge}) and by (\ref{eq:Winfty}) we have that
\begin{equation}\frac{1}{t}\int_0^tJ(s)e^{-\alpha W(s)}{\rm d}s
\end{equation}
converges almost surely to some limit and by (\ref{eq:jbar}) so does
\begin{equation}\label{eq:rhs}
\frac{\int_0^tJ(s)e^{-\alpha W(s)}{\rm d}s}{\int_0^tJ(s){\rm d}s}=\frac{\frac{1}{t}\int_0^tJ(s)e^{-\alpha W(s)}{\rm d}s}{\frac{1}{t}\int_0^tJ(s){\rm d}s}\ .
\end{equation}
Next, observe that for each $t\ge0$ (and each $\omega$ in the sample space) for which $\int_0^tJ(s){\rm d}s>0$ we have that the expression on the right hand side of (\ref{eq:rhs})
is the Laplace-Stieltjes transform of an a.s.\ nonnegative and finite random variable and thus if this ratio converges to some constant $g(\alpha)$ for each $\alpha$ then $g$ must be a Laplace-Stieltjes transform of some nonnegative (not necessarily a.s.\ finite) random variable which we denote by $W_d$. If in addition $g(\alpha)\r 1$ as $\alpha\downarrow 0$ then necessarily $g$ is the
Laplace-Stieltjes transform of a proper distribution on $\mathbb{R}_+$ and this is the case at hand as can be seen from (but is not needed for) the end result (\ref{eq:decomp}).

Finally, note that by (\ref{eq:converge}), (\ref{eq:drift}), (\ref{eq:rhs}) and the above discussion we have that
\begin{align}\label{eq:decomp1}
\varphi(\alpha)Ee^{-\alpha W(\infty)}-(\varphi(\alpha)+\eta(\alpha))p_dEe^{-\alpha W_d}\nonumber\\ =\alpha((1-p_d)\varphi'(0)-p_d\eta'(0))
\end{align}
which is equivalent to (\ref{eq:decomp}) via some obvious manipulations.\qed

\medskip
Let us now interpret (\ref{eq:decomp}). First we note that, since $\varphi'(0)>0$ then $\frac{\alpha\varphi'(0)}{\varphi(\alpha)}$ is the Laplace-Stieltjes transform of the stationary, limit and ergodic distribution associated with the process $Z_u(t)=X_u(t)+L_u(t)$ where $L_u(t)=-\inf_{0\le s\le t}X_u(s)$,
as well as the Laplace-Stieltjes transform of
the random variable $\sup_{s\ge 0}X_u(s)$. This is well known and there are quite a few proofs of this
generalized Pollaczek-Khinchin formula
in the literature, one of which is in \cite{kw1992}.

Next we observe that  from \cite{k1998}, $\frac{\eta(\alpha)}{\alpha\eta'(0)}$ is the Laplace-Stieltjes transform of the stationary excess lifetime distribution associated with the jumps of the subordinator $X_d$. For ease of reference simply observe that from
\begin{equation}
\eta(\alpha)-c_d\alpha=\int_{(0,\infty)}(1-e^{-\alpha x})\nu_d({\rm d}x)=\alpha\int_0^\infty e^{-\alpha x}\nu(x,\infty){\rm d}x
\end{equation}
and $\eta'(0)=c_d+\bar\nu_d$, where $\bar\nu_d = \int_{(0,\infty)}x\nu_d({\rm d}x)=\int_0^\infty\nu(x,\infty){\rm d}x$, we have that
\begin{equation}
\frac{\eta(\alpha)}{\alpha\eta'(0)}=\frac{c_d}{c_d+\bar\nu_d}+\frac{\bar\nu_d}{c_d+\bar\nu_d}\int_0^\infty e^{-\alpha x}\frac{\nu(x,\infty)}{\bar\nu_d}{\rm d}x
\end{equation}
which is the Laplace-Stieltjes transform of the following distribution function:
\begin{equation}
F_e(y)=\frac{c_d}{c_d+\bar\nu_d}+\frac{\bar\nu_d}{c_d+\bar\nu_d}\int_0^y \frac{\nu(x,\infty)}{\bar\nu_d}{\rm d}x
\end{equation}
for $y\ge 0$ and $F_e(y)=0$ for $y<0$. This is a somewhat generalized stationary excess lifetime distribution associated with the jumps of $X_d$.

Now introduce the random variables
$W_u,Y_e,I_l,I$:
\begin{itemize}
\item $W_u\sim\sup\{X_u(s)|\ s\ge 0\}$ with
$Ee^{-\alpha W_u}=\frac{\alpha\varphi'(0}{\varphi(\alpha)}$,
\item $Y_e\sim F_e$ with $Ee^{-\alpha Y_e}=\frac{\eta(\alpha)}{\eta'(0)\alpha}$,
\item $P(I_\ell=1)=1-P(I_\ell=0)=\pi_\ell$,
\item $P(I=1)=1-P(I=0)=\pi$,
\item $W_d$ and $W(\infty)$ are as in Theorem~\ref{th:decomp}.
\end{itemize}
Then
\begin{Theorem}\label{th:decomp1}
Under the conditions of Theorem~\ref{th:decomp}, (\ref{eq:decomp}) is equivalent to
\begin{equation}\label{eq:W_u}
W(\infty)\sim I_\ell W_u+(1-I_\ell)(I(W_u+Y_e)+W_d)\ ,
\end{equation}
where
$W_u,Y_e,I_l,I,W_d$ are assumed independent.
\end{Theorem}
We note that replacing the two instances of $W_u$ in (\ref{eq:W_u}) by two different i.i.d. random variables distributed like $W_u$ would not change the overall distribution.

One important special case of this model is when during up periods, whenever there is a positive content, the input has the same law as during down periods and the output is at a fixed rate $r>0$. That is, $\varphi(\alpha)=\alpha r-\eta(\alpha)$. A special case of this model was studied in \cite{k1998}. In this particular case it is easy to check that (as in Equation (4.12) of \cite{k1998})
\begin{equation}
1-\pi+\pi\frac{\eta(\alpha)}{\eta'(0)\alpha}\frac{\varphi'(0)\alpha}{\varphi(\alpha)}=\frac{\alpha\varphi'(0)}{\varphi(\alpha)},
\end{equation}
that is, that $I(W_u+Y_e)\sim W_u$. So in this case we have the following.
\begin{Corollary}\label{Cor:decomp2}
When $\varphi(\alpha)=\alpha r-\eta(\alpha)$ then
\begin{equation}
W(\infty)\sim W_u+(1-I_\ell)W_d,
\end{equation}
where $W_u$, $I_\ell$ and $W_d$ are independent;
and when in addition $\ell=0$ (equivalently $\tilde X(t)/t\r 0$ or $p_d=1-\pi=\frac{\varphi'(0)}{\eta'(0)+\varphi'(0)}=1-\frac{\eta'(0)}{r}$), then
\begin{equation}
W(\infty)\sim W_u+W_d\ ,
\end{equation}
where $W_u$ and $W_d$ are independent.
\end{Corollary}

We note that in Corollary~\ref{Cor:decomp2} the term $\pi=\frac{\eta'(0)}{r}$ may be referred to as the {\em traffic intensity} and is consistent with queueing theory.

\begin{Remark}\rm
Throughout this and the following section we are focussing on almost sure convergence.
However, throughout, most ``almost sure'' statements could be trivially replaced by ``in probability'' without changing anything else (simply by looking at subsequences that converge a.s.). We are not aware of related applications where the convergence is in probability but not almost surely and thus did not see a point in making this issue more precise.
\end{Remark}
\begin{Remark}\rm
In \cite{kw1991} the focus is on convergence in distribution rather than long run a.s. convergence. As in the previous remark, we could follow the same ideas with similar proofs (but with more restrictive assumptions). We chose to leave this out as, given what follows, and what is already available in \cite{kw1991}, it may be considered an exercise.
\end{Remark}

\section{How to interpret $W_d$?}\label{sec:Wd}
In this section we identify the non-standard component in the decomposition
of Theorem~\ref{th:decomp}, associated with down periods.
In particular, we will express the Laplace-Stieltjes transform
$\frac{\eta(\alpha)}{\alpha \eta'(0)} Ee^{-\alpha W_d}$ of $Y_e + W_d$ in terms of
the transforms of the workloads at the ends of up and down periods.

We recall that under the assumptions of Theorem~\ref{th:decomp},
\begin{equation}
Ee^{-\alpha W_d}=\lim_{t\r\infty}\frac{\int_0^te^{-\alpha W(s)}J(s){\rm d}s}{\int_0^tJ(s){\rm d}s}
\end{equation}
and since for every nonnegative random variable $V$ we have that $e^{-\alpha V}=\alpha \int_0^\infty e^{-\alpha x}\one{V\le x}{\rm d}x$, then also here
\begin{equation}
\frac{\int_0^te^{-\alpha W(s)}J(s){\rm d}s}{\int_0^tJ(s){\rm d}s}=\alpha \int_0^\infty e^{-\alpha x}\frac{\int_0^t\one{W(s)\le x}J(s){\rm d}s}{\int_0^tJ(s){\rm d}s}{\rm d}x
\end{equation}
and thus, a.s., $\frac{\int_0^t\one{W(s)\in \cdot}J(s){\rm d}s}{\int_0^tJ(s){\rm d}s}$ (probability distribution valued process) converges in distribution to $W_d$.
This holds in particular if we replace $t$ by $S_n$. In this case $\int_0^{S_n}J(s){\rm d}s=\sum_{k=1}^n X_k$ and thus we have that
\begin{equation}
\frac{\int_0^{S_n}\one{W(s)\in\cdot}J(s){\rm d}s}{\int_0^{S_n}J(s){\rm d}s}=\frac{\sum_{k=1}^n\int_0^{X_k}\one{W(T_{k-1}+s)\in \cdot}{\rm d}s}{\sum_{k=1}^nX_k}
\end{equation}
where for $s\in[0,X_n)$ we have that
\begin{equation}
W(T_{n-1}+s)=W(T_{n-1})+X_d(T_{n-1}+s)-X_d(T_{n-1})\
\end{equation}
and thus
\begin{equation}
\int_0^{X_n}e^{-\alpha W(T_{n-1}+s)}{\rm d}s=e^{-\alpha W(T_{n-1})}\int_0^{X_n}e^{-\alpha(X_d(T_{n-1}+s)-X_d(T_{n-1}))}{\rm d}s\ .
\end{equation}
Now, since $T_{n-1},S_n$ are stopping times with respect to $\left\{\mathcal{F}_t|\ t\ge0\right\}$, $X_n$ is a stopping time with respect to
$\left\{\mathcal{F}_{T_{n-1}+t}|\ t\ge0\right\}$ (of course, not with respect to the original filtration in general). Moreover, $W(T_{n-1})\in \mathcal{F}_{T_{n-1}}$ and by the strong Markov property $X_d^{T_{n-1}}\equiv\{X_d(T_{n-1}+t)-X_d(T_{n-1})|\ t\ge 0\}$ is a subordinator with respect to $\left\{\mathcal{F}_{T_{n-1}+t}|\ t\ge0\right\}$ with exponent $\eta$ (that is, distributed like $X_d$) and is independent of $F_{T_{n-1}}$ (thus, of $W(T_{n-1})$). Thus from \cite{kw1992} we have that
\begin{equation}
-\eta(\alpha)\int_0^te^{-\alpha X_d^{T_{n-1}}(s)}{\rm d}s+1-e^{-\alpha X_d^{T_{n-1}}(t)}
\end{equation}
is a zero mean martingale with respect to $\left\{\mathcal{F}_{T_{n-1}+t}|\ t\ge0\right\}$ and thus by the optional stopping theorem together with monotone and bounded convergence where appropriate we have with
\begin{equation}
\Delta_n=-\eta(\alpha)\int_0^{X_n}e^{-\alpha X_d^{T_{n-1}}(s)}{\rm d}s+1-e^{-\alpha X_d^{T_{n-1}}(X_n)} ,
\end{equation}
that
$E[\Delta_n|\mathcal{F}_{T_{n-1}}]=0$. Moreover, from Lemma 3 of \cite{KB} and the fact that $M(t)^2-[M,M](t)$ is a (zero mean) martingale, we can conclude that when $X_n$ is a.s. finite then
\begin{equation}
E[\Delta_n^2|\mathcal{F}_{T_{n-1}}]=(2\eta(\alpha)-\eta(2\alpha))E\left[\left.\int_0^{X_n} e^{-2\alpha X_d^{T_{n-1}}(s)}{\rm d}s\right|\mathcal{F}_{T_{n-1}}\right]
\end{equation}
and in the same way that led to $E[\Delta_n|\mathcal{F}_{T_{n-1}}]=0$, by substituting $2\alpha$ instead of $\alpha$, we have that
\begin{equation}
\eta(2\alpha)E\left[\left.\int_0^{X_n}e^{-2\alpha X_d^{T_{n-1}}(s)}{\rm d}s\right|\mathcal{F}_{T_{n-1}}\right]=1-E\left[\left.e^{-2\alpha X_d^{T_{n-1}}(X_n)}\right|\mathcal{F}_{T_{n-1}}\right]
\end{equation}
and we conclude that
\begin{equation}
E[\Delta_n^2|\mathcal{F}_{T_{n-1}}]=\left(\frac{2\eta(\alpha)}{\eta(2\alpha)}-1\right)\left(1-E\left[\left.e^{-2\alpha X_d^{T_{n-1}}(X_n)}\right|\mathcal{F}_{T_{n-1}}\right]\right)\ .
\end{equation}
In particular, upon multiplying by $e^{-\alpha W(T_{n-1})}\in\mathcal{F}_{T_{n-1}}$, we have that
\begin{equation}
\sum_{k=1}^n e^{-\alpha W(T_{k-1})} \Delta_k
\end{equation}
is a zero mean martingale, where
\begin{equation}
E\left[\left.\left(e^{-\alpha W(T_{k-1})} \Delta_k\right)^2\right|\mathcal{F}_{T_{k-1}}\right]\le \frac{2\eta(\alpha)}{\eta(2\alpha)}-1<\infty\ .
\end{equation}
 It is well known (cf. Theorem 3 on p.\ 243 of \cite{Feller}) that an $L^2$ martingale $M_n$ satisfying
 \begin{equation}
 \sum_{k=1}^\infty \frac{E(M_k- M_{k-1})^2}{k^2}<\infty
 \end{equation}
 also satisfies $M_n/n\r0$ a.s. and in $L^2$ and thus
\begin{equation}
\frac{1}{n}\sum_{k=1}^n e^{-\alpha W(T_{k-1})} \Delta_k\r0
\end{equation}
a.s. and in $L^2$ and we finally have the following.
\begin{Theorem}
Under the assumptions of Theorem~\ref{th:decomp},
\begin{equation}
\frac{1}{n}\left(-\eta(\alpha)\int_0^{S_n} e^{-\alpha W(s)}J(s){\rm d}s+\sum_{k=1}^n\left(e^{-\alpha W(T_{k-1})}-e^{-\alpha W(S_k)}\right)\right)\r 0
\end{equation}
a.s. and in $L^2$, and if in addition $p_d>0$ then
\begin{equation}
\frac{-\eta(\alpha)\int_0^{S_n} e^{-\alpha W(s)}J(s){\rm d}s+\sum_{k=1}^n\left(e^{-\alpha W(T_{k-1})}-e^{-\alpha W(S_k)}\right)}{\int_0^{S_n}J(s){\rm d}s}\r 0
\end{equation}
and thus
\begin{equation}
\frac{\sum_{k=1}^n\left(e^{-\alpha W(T_{k-1})}-e^{-\alpha W(S_k)}\right)}{\sum_{k=1}^nX_k}\r \eta(\alpha)Ee^{-\alpha W_d} .
\end{equation}
\end{Theorem}
Now, note that from $\frac{1}{t}\int_0^tJ(s){\rm d}s\r p_d>0$, if also $T_n/n\r \mu>0$ a.s. (and thus also $S_n/n\r \mu$) then
\begin{equation}
\frac{1}{n}\sum_{k=1}^nX_k=\frac{S_n}{n}\frac{1}{S_n}\int_0^{S_n}\one{J(s)}{\rm d}s\r\mu p_d>0
\end{equation}
and thus
\begin{equation}
\frac{1}{n}\sum_{k=1}^n\left(e^{-\alpha W(T_{k-1})}-e^{-\alpha W(S_k)}\right)
\end{equation}
converges a.s. In particular, we have:
\begin{Theorem}
Under the assumptions of Theorem~\ref{th:decomp}, if $p_d>0$ and $T_n/n\r\mu>0$ a.s., then
$\frac{1}{n}\sum_{k=1}^ne^{-\alpha W(S_k)}\r Ee^{-\alpha W_+}$ a.s. for some nonnegative random variable $W_+$ if and only if $\frac{1}{n}\sum_{k=1}^ne^{-\alpha W(T_{k-1})}\r Ee^{-\alpha W_-}$ a.s. for some nonnegative random variable $W_-$ and we have that
\begin{equation}
\frac{Ee^{-\alpha W_-}-Ee^{-\alpha W_+}}{\alpha \eta'(0)\mu p_d}=\frac{\eta(\alpha)}{\alpha\eta'(0)}Ee^{-\alpha W_d}\ .
\end{equation}
Moreover if any two of $EW_d,EW_-,EW_+$ are finite, then so is the third and we have that
\begin{equation}\label{eq:pm}
\frac{Ee^{-\alpha W_-}-Ee^{-\alpha W_+}}{\alpha (EW_+-EW_-)}=\frac{\eta(\alpha)}{\alpha\eta'(0)}Ee^{-\alpha W_d}\ .
\end{equation}
\end{Theorem}
The above theorem gives more insight into the distribution and meaning of
$W_d + Y_e$, by relating this sum to the random variables $W_+$ and $W_-$ which successively represent the workload at the ends of down and up periods.
 For more details regarding the left side of (\ref{eq:pm}), see Theorems~5.1 and 5.2 of \cite{kw1991}. In particular, it is a Laplace-Stieltjes transform of a bona fide distribution if and only if $W_-$ is stochastically smaller than $W_+$. This form was also observed and discussed in the M/G/1 queue setting in \cite{s1988}. Finally, if there are enough assumptions to assure that $W_-$ and $U=W_+-W_-$ are independent then the left side of (\ref{eq:pm}) becomes
 \begin{equation}
 Ee^{-\alpha W_-}\ \frac{1-Ee^{-\alpha U}}{\alpha EU}\ .
 \end{equation}
That is, it is the transform of a sum of two independent random variables, the first is $W_-$ and the second has the stationary residual lifetime distribution of $U$. If we denote this variable by $U_e$ then we have the following decomposition
\begin{equation}
W_-+U_e\sim W_d+Y_e\ ,
\end{equation}
where we recall that $Y_e$ has the transform $\frac{\eta(\alpha)}{\alpha\eta'(0)}$ and the variables on either side are assumed independent. The special case where this kind of independence (between $W_-$ and $U$) occurs is discussed in the M/G/1 queue setting in \cite{fc1985}. We also refer the reader to Theorem~4.1 and its proof in \cite{k1998} for the special case considered there.

We recall that by Theorem~\ref{th:decomp1},
\begin{equation}
W(\infty)\sim I_\ell W_u+(1-I_\ell)(I(W_u+Y_e)+W_d)\ .
\end{equation}
Thus, replacing $Y_e$ by an independent $Y_e^1\sim Y_e$ and adding $Y_e$ on both sides we have that
\begin{equation}
W(\infty)+Y_e\sim I_\ell (W_u+Y_e)+(1-I_\ell)(I(W_u+Y_e^1)+(W_d+Y_e))\ .
\end{equation}
With $W_\pm\sim W_d+Y_e$ (a random variable with LST given by the left side of (\ref{eq:pm})) this implies that
\begin{equation}
W(\infty)+Y_e\sim I_\ell (W_u+Y_e)+(1-I_\ell)(I(W_u+Y_e^1)+W_\pm)\ ,
\end{equation}
where the expressions on either side of the equation are independent. Finally, replacing $Y^1_e$ on the right by $Y_e$ does not change the distribution (due to the indicator $I_\ell$) so that
\begin{equation}
W(\infty)+Y_e\sim I_\ell(W_u+Y_e)+(1-I_\ell)(I(W_u+Y_e)+W_\pm)\ ,
\end{equation}
where again all variables appearing on the expressions on either side of the equation are assumed independent so that only their marginal distribution matters.
In the special case of $\varphi(\alpha)=\alpha r-\eta(\alpha)$ we
replace $I(W_u+Y_e)$ on the right by $W_u$ (see Corollary~\ref{Cor:decomp2}) and obtain
\begin{eqnarray}
W(\infty)+Y_e&\sim &I_\ell(W_u+Y_e)+(1-I_\ell)(W_u+W_\pm)\nonumber \\ &=&W_u+I_\ell Y_e+(1-I_\ell)W_\pm\ ,
\end{eqnarray}
and in particular when $\pi_\ell=0$
\begin{equation}
W(\infty)+Y_e\sim W_u+W_\pm\ ,
\end{equation}
where, again, throughout all random variables appearing in the expressions on either sides of the equations are assumed independent.

\section{Applications to polling systems}\label{sec:polling}
In this section we relate our decomposition results to decomposition results for so-called polling systems.
A polling system is a single-server multi-queue system, in which the server visits the queues one at a time, typically in a cyclic order.
The service discipline at each queue specifies the duration of a visit.
E.g., under the exhaustive service discipline, the server visits a queue until it has become empty;
under the $1$-limited discipline, it serves exactly one customer during a visit.
In many applications (e.g., in production systems, where the server is a machine and the customers of a queue are orders
of a particular type) it is natural to have nonnegligible switchover times from one queue to the next.
Stimulated by a wide variety of applications (not only production systems,  but also computer- and communication systems,
traffic lights, repair systems), polling models have been extensively studied.
It is almost always assumed that the input processes to the queues are independent Poisson processes.
For such a situation, it was proven in \cite{BG} that the steady state total workload in the polling system
{\em with} switchover times
can be decomposed into
two independent quantities, viz. (i) the workload in the corresponding polling system {\em without} switchover times,
and (ii)  the steady state total amount of work at an epoch the server is not working.
Item (i) is the workload in an $M/G/1$ queueing system; the distribution of item (ii) was determined
for a few service disciplines in \cite{BTT}.
In \cite{BKK} the joint steady state workload distribution at arbitrary epochs was expressed
in the joint queue length distribution at visit beginning and visit completion epochs.
The latter distributions are known for
certain polling models, in particular, for polling models in which
the service discipline at all queues is of so-called branching type.

The cyclic polling model of \cite{BG} was generalized in
\cite{BGW} to the case of a fixed non-cyclic visit order of the queues; again a work decomposition result was derived.
A further generalization is contained in
\cite{boxma}. That paper considers a single-server multi-class system
with a work-conserving scheduling discipline as long as the server {\em is} serving and with
a service interruption process (which could correspond to switchover times in a polling system)
that does not affect the amount of service time given to a customer
or the arrival time of any customer. Furthermore, the arrival process is a batch Poisson process
that allows correlations between
the numbers of arrivals of the various customer types in a batch.
Again a decomposition result was proven: the steady state workload in the model with interruptions
is in distribution equal to the sum of two independent quantities,
viz. (i) the steady state workload in the corresponding model without interruptions,
and (ii) the steady state amount of work at an epoch
in which the server is not serving.

Another extension of the cyclic polling model with independent Poisson arrivals was recently studied in
\cite{BIKM}. It considers a cyclic polling system with $N$ queues, extending the Poisson arrival process to an $N$-dimensional L\'evy subordinator (so the sample paths are non-decreasing in all coordinates).
If a particular queue is being served, then the workload level at that queue behaves as a spectrally positive
L\'evy process with a negative drift. Another special feature of the model is that the L\'evy input process
changes at polling and switching instants. A restrictive assumption is that the service discipline at each queue
is of branching type. That assumption implies that the $N$-dimensional workload process
at successive instants that the server arrives at the first queue is a {\em Jirina} process,
which is a multi-type continuous-state branching process.
The joint steady state workload distribution at such epochs, and subsequently also at arbitrary epochs,
is determined in \cite{BIKM}; no workload decomposition is derived.
A special case (constant fluid input at all queues) had been studied by Czerniak and Yechiali \cite{CY}, who also obtained the joint workload
distribution at arbitrary epochs. In Section 4 of their paper
they point out that, if there is a workload decomposition, the term "(i)" without switchover times is zero because the outflow
is larger than the inflow during visit times.

In Section~\ref{sec:4} of the present paper, we derive workload transforms and workload decompositions in a system
that alternates between up and down periods.  The input process is one  L\'evy process
$X_u$ during up periods and another L\'evy process $X_d$ during down periods.
Our Theorem~\ref{th:decomp} generalizes exact workload transform results in \cite{BKK} and \cite{BTT},
where the input process is a sum of independent compound Poisson processes, to the case
of a L\'evy input process.
It complements the exact workload transform result of \cite{BIKM} in the sense that it
only gives total workload and
does not give a joint transform, but that it does allow more general visit disciplines.
Our assumption on the {\em up} and {\em down} periods (visit times and interruptions), viz., the assumption that
$0=T_0 \leq S_1 \leq T_1 \leq S_2 \leq T_2 \dots$ is an increasing sequence of a.s. finite stopping times,
in particular includes non-branching service disciplines.
Our Corollary~\ref{Cor:decomp2} generalizes/complements decomposition results for total workload in \cite{BG,boxma}
for polling systems and, more generally, single-server multi-class systems with interruptions.
Our L\'evy input process generalizes the (batch) Poisson processes of those and other polling papers.

In fact, due to our general setup it seems that under appropriate stability conditions, decomposition results would hold for quite general polling mechanisms. Some examples are cases were the lengths of the switching times depend on the state of the system in various ways (e.g., shorter switching when certain queues are large), or when the decision of when to leave a certain queue may depend on the overall information of the system rather than following a fixed mechanism.

\end{document}